\documentclass[10pt,a4paper]{article}
\usepackage{graphics, graphicx}     

\usepackage{cite}
\usepackage{amsmath}
\usepackage{amssymb}
\usepackage{amsthm}
\theoremstyle{plain}
\newtheorem{Theorem}{Theorem}[section] %
\newtheorem{Lemma}{Lemma}[section]

\theoremstyle{definition}
\newtheorem{Remark}{Remark}[section]

\theoremstyle{definition}

\newtheorem{Problem}{Problem}[section]

{\par\noindent{\it Proof of}} 
{\hfill$\vspace{5mm}\scriptstyle\blacksquare$} 

\numberwithin{equation}{section} 
\numberwithin{figure}{section} 
\numberwithin{table}{section} 

\begin{document}

\setcounter{page}{1}

\markboth{M.I. Isaev}{Energy and regularity dependent stability estimates for near-field  inverse scattering in multidimensions}

\title{Energy and regularity dependent stability estimates for  near-field inverse scattering in multidimensions}
\date{}
\author{ {\bf M.I. Isaev} \\ Centre de Math\'ematiques Appliqu\'ees, Ecole Polytechnique,\\
91128 Palaiseau, France\\
e-mail: \tt{isaev.m.i@gmail.com}}

\maketitle
\begin{abstract}
We prove new global H\"older-logarithmic stability estimates for the near-field  inverse scattering  problem in dimension $d\geq 3$.
Our estimates
are given in uniform norm for coefficient difference and related
stability efficiently increases with increasing energy and/or coefficient regularity. 
In addition, a global logarithmic stability estimate for this inverse problem in dimension $d=2$ is also given.
\end{abstract}

\section{Introduction}
We consider the  Schr\"odinger equation
\begin{equation}\label{eq} 
	L\psi = E\psi, \ \  \  L =-\Delta   + v(x), \ \ \  x \in  \mathbb{R}^d, \ \  d\geq 2,
\end{equation}
where
\begin{equation}\label{eq_c}
 \begin{aligned}
 	 &v \text{ is real-valued},  \ \ v\in \mathbb{L}^{\infty}(\mathbb{R}^d),\\
  	&v(x) = O(|x|^{-d-\varepsilon}),  \ \ |x|\rightarrow \infty, \ \   \text{for some } \varepsilon>0.
 \end{aligned}
\end{equation}
We consider the resolvent $R(E)$ of the  Schr\"odinger operator $L$ in $\mathbb{L}^2(\mathbb{R}^d)$:
\begin{equation}
	R(E) = (L-E)^{-1},\ \ \ \ E\in \mathbb{C} \setminus \sigma(L),
\end{equation}
where $\sigma(L)$ is the spectrum of $L$  in $\mathbb{L}^2(\mathbb{R}^d)$.
We assume that $R(x,y,E)$ denotes the Schwartz kernel of $R(E)$ as of an integral operator. 
We consider also 
\begin{equation}
	R^+(x,y,E) = R(x,y,E + i0), \ \ \ x,y \in \mathbb{R}^d, \ \  E\in \mathbb{R}_+.
\end{equation}
We recall that 
in the framework of equation \eqref{eq}
the function $R^+(x,y,E)$ 
 describes scattering of the spherical waves  
 \begin{equation}
 R^+_0(x,y,E) = -\frac{i}{4} \left(\frac{\sqrt{E}}{2\pi |x-y|}\right)^{\frac{d-2}{2}}H^{(1)}_{\frac{d-2}{2}}(\sqrt{E}|x-y|),
 \end{equation}
 generated by a source at $y$ (where $H^{(1)}_\mu$ is the Hankel function of the first kind of order $\mu$).
We recall also that $R^+(x,y,E)$ is the Green function for $L-E$, $E\in\mathbb{R}_+$, with  the Sommerfeld radiation condition at infinity. 

In addition, the function
\begin{equation}
	\begin{aligned}
		S^+(x,y,E) = R^+(x,y,E) - R^+_{0}(x,y,E),\\
			  x, y \in \partial B_r,\  E\in\mathbb{R}_+, \ r\in\mathbb{R}_+,
	\end{aligned}
\end{equation}
is considered as near-field scattering data for equation \eqref{eq},
 where  $B_r$ is the open ball of radius $r$ centered at $0$.

We consider, in particular, the following near-field inverse scattering problem for equation (\ref{eq}): 

\begin{Problem}
 Given $S^+$ on $\partial B_r\times \partial B_r$ for some fixed $r,E\in\mathbb{R}_+$, find $v$ on $B_r$. 
\end{Problem}

 This problem can be considered under the assumption that $v$ is a priori known on $\mathbb{R}^d\setminus B_r$.
Actually, in the present paper  we consider Problem 1.1 under the assumption that $v\equiv 0$ on $\mathbb{R}^d\setminus B_r$ 
for some fixed $r\in \mathbb{R}_+$. Below in this paper we always assume that this additional condition is fulfilled.

It is well-known that the near-field scattering data of Problem 1.1 uniquely and efficiently determine  the scattering 
amplitude $f$ for equation \eqref{eq} at fixed  energy $E$,  see \cite{Berezanskii1958}. Therefore, approaches of \cite{ABR2008}, \cite{Buckhgeim2008},
\cite{BAR2009}, \cite{ER1995}, \cite{Henkin1987}, \cite{IsaevFunc}, 
 \cite{Novikov1988}, \cite{Novikov1994}, \cite{Novikov1999}, \cite{Novikov2005+}, \cite{Stefanov1990}  can be applied to Problem 1.1 via this reduction.

In addition,  it is also known  that the near-field data of Problem 1.1  uniquely
determine the Dirichlet-to-Neumann map in the case when $E$ is not a Dirichlet eigenvalue for operator $L$ in $B_r$, 
see \cite{Nachman1988}, \cite{Novikov1988}. Therefore, approaches
of \cite{Alessandrini1988}, \cite{Buckhgeim2008}, \cite{IN2012++}, \cite{Isakov2011}, \cite{Mandache2001}, \cite{Novikov1988}, \cite{Novikov2005}-\cite{NS2012}, \cite{SU1987}
can be also applied to Problem 1.1 via this reduction.

However, in some case it is much more optimal to deal with Problem 1.1 directly, see, for example, logarithmic stability results of \cite{HH2001} for Problem 1.1 in dimension $d = 3$. A principal improvement of estimates of \cite{HH2001} was given recently in \cite{IN2012+++}:
stability of \cite{IN2012+++} efficiently increases with increasing regularity of $v$.

Problem 1.1 can be also considered as an example of ill-posed problem: 
see \cite{LR1986}, \cite{BK2012} for an introduction to this theory.

In the present paper we continue studies of \cite{HH2001}, \cite{IN2012+++}. 
We give new global H\"older-logarithmic stability estimates
for Problem 1.1 in dimension $d\geq 3$, see Theorem 2.1. Our estimates are given in uniform norm for coefficient
difference and related stability efficiently increases with increasing energy and/or
coefficient regularity. Results of such a type for the Gel'fand inverse problem  were obtained recently in \cite{IN2012++} for $d\geq 3$ and in \cite{S2012} for $d=2$. 

In addition, we give also global logarithmic stability estimates
for Problem 1.1 in dimension $d=2$, see Theorem 2.2.


	
\section{Stability estimates}

We recall that if  $v$ satisfies  \eqref{eq_c} and $\mbox{supp}\, v \subset B_{r_1}$ for some $r_1>0$, then
 \begin{equation}
	S^+(E) \text { is bounded  in } \mathbb{L}^{2}(\partial B_r \times \partial B_r) \text{ for any $r>r_1$,}
\end{equation}
  where  $S^+(E)$ is the near-field scattering data of  $v$ for equation \eqref{eq}
  with $E>0$, for more details see, for example, Section 2  of \cite{HH2001}.

\vspace{1mm}

\noindent
{\it 2.1 Estimates for $d \geq 3$}

\vspace{1mm}

In this subsection we assume for simplicity that
\begin{equation}\label{assumption}
	\begin{aligned}
		 v \in \mathbb{W}^{m,1}(\mathbb{R}^d) \text{ for some } m > d, \\ 
		  		 	v \text{ is real-valued},\\
		 \mbox{supp}\, v \subset B_{r_1} \  \text{ for some $r_1>0$},
	\end{aligned}
\end{equation}
where 
\begin{equation}\label{2.5}
	\mathbb{W}^{m,1}(\mathbb{R}^d) = \{v:\  \partial^J v \in \mathbb{L}^1(\mathbb{R}^d),\  |J| \leq m \},\ m\in \mathbb{N}\cup 0,
\end{equation}
where 
\begin{equation}
J \in (\mathbb{N}\cup 0)^d,\ |J| = \sum\limits_{i=1}\limits^{d}J_i,\ \partial^J v(x) 
= \frac{\partial^{|J|} v(x)}{\partial x_1^{J_1}\ldots \partial x_d^{J_d}}.
\end{equation}
Let
\begin{equation}
	||v||_{m,1} = \max\limits_{|J|\leq m} ||\partial^J v||_{\mathbb{L}^1(\mathbb{R}^d)}.
\end{equation}

Note that (\ref{assumption}) $\Rightarrow$ (\ref{eq_c}).

\begin{Theorem}\label{Theorem_2.1} Let $E>0$ and $r>r_1$ be given constants.  
Let  dimension $d\geq 3$ and potentials $v_1$, $v_2$ satisfy  \eqref{assumption}.
Let $||v_{j}||_{m,1} \leq N,\  j = 1,2,$ for some $N>0$. 
Let $S_1^+(E)$ and $S_2^+(E)$ denote 
the near-field scattering data for 
$v_1$ and  $v_2$, respectively. Then for  
  $\tau \in(0,1)$ and any $s \in [0, s^*]$ the following estimate holds: 
\begin{equation}\label{eq_t2}
	||v_2 - v_1||_{L^\infty(B_r)} \leq C_1 (1+E)^{\frac{5}{2}} \delta^\tau + 
	C_2 (1+E)^{\frac{s- s^*}{2}}\left(\ln\left(3+\delta ^{-1}\right)\right)^{-s}, 
\end{equation}
where  $s^* = \frac{m-d}{d}$, $\delta  = ||S_1^+(E) - S_2^+(E)||_{\mathbb{L}^{2}(\partial B_r \times \partial B_r)}$,  and constants 
$C_1,C_2>0$ depend only on $N$,  $m$, $d$, $r$,  $\tau$.
\end{Theorem}

Proof of Theorem \ref{Theorem_2.1} is given in Section 5. This proof is based on results presented in Sections 3, 4.

\vspace{1mm}

\noindent
{\it 2.2 Estimates for $d = 2$}

\vspace{1mm}

In this subsection we assume for simplicity that
\begin{equation}\label{assumption2}
	\begin{aligned}
		 v \text{ is real-valued}, \ \ v \in C^{2}(\overline{B}_{r_1}),  \\ \mbox{supp}\, v \subset B_{r_1} \  \text{ for some $r_1>0$}.
		 \end{aligned}
\end{equation}
Note also that (\ref{assumption2}) $\Rightarrow$ (\ref{eq_c}).

\begin{Theorem}\label{Theorem_2.2}
Let $E>0$ and $r>r_1$ be given constants.  
Let  dimension $d = 2$ and potentials $v_1$, $v_2$ satisfy  \eqref{assumption2}.
Let $||v_{j}||_{C^2(B_{r})} \leq N,\  j = 1,2,$ for some $N>0$. 
Let $S_1^+(E)$ and $S_2^+(E)$ denote 
the near-field scattering data for 
$v_1$ and  $v_2$, respectively. Then 
\begin{equation}\label{estimation_2}
	\begin{aligned}
	||v_1 - v_2||_{L^\infty(B_r)} \leq 
	C_3 \left(\ln\left(3+\delta^{-1}\right)\right)^{-3/4}
	\left( \ln \left(3\ln\left(3+\delta^{-1}\right)\right) \right)^2,
	\end{aligned}
\end{equation}
where  $\delta  = ||S_1^+(E) - S_2^+(E)||_{\mathbb{L}^{2}(\partial B_r \times \partial B_r)}$  and constant 
$C_3>0$ depends only on $N$, $m$, $r$.
\end{Theorem}

Proof of Theorem \ref{Theorem_2.2} is given in Section 7. This proof is based on results presented in Sections 3, 6.

\vspace{1mm}

\noindent
{\it 2.3 Concluding remarks}
\begin{Remark}
	The logarithmic stability estimates for Problem 1.1 of \cite{HH2001} and \cite{IN2012+++} follow from estimate \eqref{eq_t2}
	for $d=3$ and $s=s^*$.
	Apparently, using the methods of \cite{Novikov2009}, \cite{Novikov2011} it is possible to improve estimate \eqref{eq_t2} 
for $s^* = m-d$.
\end{Remark}

\begin{Remark}
	In the same way as in \cite{HH2001} and \cite{IN2012+++} for dimesnsion $d=3$,	
	using estimates \eqref{eq_t2} and \eqref{estimation_2}, one can obtain logarithmic stability estimates for the reconstruction 
	of a potential $v$ from the inverse scattering amplitude $f$ for any $d\geq 2$.
\end{Remark}

\begin{Remark}
	Actually, in the proof of Theorem 2.1 we obtain the following
estimate (see formula \eqref{5.18}):
\begin{equation}\label{Remark_eq}
\|v_1 - v_2\|_{\mathbb{L}^{\infty}(B_r)} \leq C_4 (1+E)^{2}\sqrt{E+\rho^2}\, e^{2\rho (r+1)} \delta+  
 	C_5(E+\rho^2)^{-\frac{m-d}{2d}},
\end{equation}
where constants $C_4, C_5>0$ depend only on $N$,  $m$, $d$, $r$ and the parameter $\rho>0$
is such that $E+\rho^2$ is sufficiently large: $E+\rho^2\geq C_6(N,r,m)$.  Estimate of Theorem 2.1 follows from estimate \eqref{Remark_eq}.
\end{Remark}

\section{Alessandrini-type identity for near-field scattering}
In this section we always assume that assumptions of Theorems \ref{Theorem_2.1} and \ref{Theorem_2.2}
are fulfilled (in the cases of dimension $d\geq 3$ and $d=2$, respectively). 

Consider the operators $\hat{\mbox{R}}_{j}$, $j=1,2$, defined as follows
\begin{equation}
	(\hat{\mbox{R}}_{j} \phi) (x) = \int\limits_{\partial B_r} R_j^+(x,y,E)  \phi(y) dy, \ \ \ x\in \partial B_r, \  \ j=1,2.
\end{equation}
Note that
\begin{equation}\label{eq4.2}
	\|\hat{\mbox{R}}_{1} - \hat{\mbox{R}}_{2} \|_{\mathbb{L}^2(\partial B_r)} \leq 
	\|S_1^+(E)  - S_2^+(E)   \|_{\mathbb{L}^2(\partial B_r) \times \mathbb{L}^2(\partial B_r)}.
\end{equation}
We recall that (see \cite{HH2001}) for any  functions  $\phi_1,\phi_2\in C(\mathbb{R}^d)$, 
 sufficiently regular in  $\mathbb{R}^d\setminus \partial B_r $ and satisfying
\begin{equation}\label{4.3}
	\begin{aligned}
	-\Delta\phi   + v(x) \phi = E\phi, \ \ \text{ in } \mathbb{R}^d \setminus \partial B_r, \\
	\lim\limits_{|x|\rightarrow +\infty} |x|^{\frac{d-1}{2}} \left(\frac{\partial}{\partial|x|}\phi - i\sqrt{E}\phi\right)=0,\\
	\end{aligned}
\end{equation} 
with $v=v_1$ and $v=v_2$, respectively,
the following identity holds:
\begin{equation}\label{Id}
	\begin{aligned}
	\int\limits_{B_r}(v_2 -v_1) \phi_1 \phi_2 dx &=\\ =\int \limits_{\partial B_r} &\left(\frac{\partial{\phi_1}}{\partial \nu_+} - 
	\frac{\partial{\phi_1}}{\partial \nu_-}\right)
	\left[
	\left(\hat{\mbox{R}}_{1} - \hat{\mbox{R}}_{2}\right) 
	\left(\frac{\partial{\phi_2}}{\partial \nu_+} - \frac{\partial{\phi_2}}{\partial \nu_-}\right)
	\right]dx,
	\end{aligned}
\end{equation}
where where $\nu_{+}$ and $\nu_{-}$  are the
outward  and inward normals to $\partial B_r$, respectively.

\begin{Remark} 
The identity \eqref{Id} is similar to the Alessandrini identity (see Lemma 1 of [1]),  where the Dirichlet-to-Neumann maps are considered  instead of operators $\hat{\mbox{R}}_j$. 
\end{Remark}

To apply identity \eqref{Id} to our considerations, we use also the following lemma:
\begin{Lemma}\label{Lemma_4.1}
	 Let $E,r>0$ and  $d\geq 2$.  Then, there is 
	a positive constant $C_7$ (depending only on $r$ and $d$) such that for any  $\phi\in C(\mathbb{R}^d \setminus B_r)$ satisfying
	\begin{equation}
	\begin{aligned}
	-\Delta \phi   = E\phi, \ \ \text{ in } \mathbb{R}^d \setminus \overline{B}_r, \\
	\lim\limits_{|x|\rightarrow +\infty} |x|^{\frac{d-1}{2}} \left(\frac{\partial}{\partial|x|}\phi - i\sqrt{E}\phi\right)=0,\\
	\phi|_{\partial B_r} \in \mathbb{H}^1(\partial B_r),
	\end{aligned}
\end{equation}
the following inequality holds:
	\begin{equation}\label{eq4.3}
		\left\|\left. \frac{\partial\phi}{\partial\nu_+}\right|_{\partial B_r}\right\|_{\mathbb{L}^2(\partial B_r)} \leq C_7 (1+E)
		\left\| \phi|_{\partial B_r} \right\|_{\mathbb{H}^1(\partial B_r)},
	\end{equation}
where $\mathbb{H}^1(\partial B_r)$ denotes the standart Sobolev space on $\partial B_r$.
\end{Lemma}
The proof of Lemma \ref{Lemma_4.1} is given in Section 8.

\section{Faddeev functions}
In dimension $d\geq 3$,
we consider the Faddeev functions $h$, $\psi$, $G$ (see \cite{Faddeev1965}, \cite{Faddeev1974}, \cite{Henkin1987}, \cite{Novikov 1988}): 
\begin{equation}\label{4.4}
	h(k,l) = (2\pi)^{-d} \int\limits_{\mathbb{R}^d} e^{-ilx}v(x) \psi(x,k) dx,
\end{equation}
where $k,l \in \mathbb{C}^d, \ k^2=l^2,\ \mbox{Im}\, k = \mbox{Im}\, l \neq 0$,
\begin{equation}\label{4.1}
	\psi(x,k) = e^{ikx} + \int\limits_{\mathbb{R}^d} G(x-y,k) v(y)\psi(y,k) dy, 
\end{equation}
\begin{equation}\label{4.2}
	G(x,k) = e^{ikx} g(x,k), \ \ \ g(x,k) = - (2\pi)^{-d} \int\limits_{\mathbb{R}^d} \frac{e^{i\xi x} d\xi}{\xi^2 + 2k\xi},
\end{equation}
where $x\in \mathbb{R}^d$, $k\in \mathbb{C}^d$, $\mbox{Im}\, k\neq 0$, $d\geq 3$,
	 
One can consider  (\ref{4.4}), (\ref{4.1}) assuming that 
\begin{equation}\label{4.6}
	\begin{aligned}
	v \text{ is a sufficiently regular function on } \mathbb{R}^d \\
	\text{ with suffucient decay at infinity.}
	\end{aligned}
\end{equation}
For example, in connection with Theorem \ref{Theorem_2.1},  we consider  (\ref{4.4}), (\ref{4.1}) assuming that
\begin{equation}\label{4.7}
	v \in \mathbb{L}^{\infty}(B_r), \ \ \ v \equiv 0 \text{ on } \mathbb{R}\setminus B_r.
\end{equation}

We recall that (see \cite{Faddeev1965}, \cite{Faddeev1974}, \cite{Henkin1987}, \cite{Novikov 1988}): 
\begin{equation}\label{4.8}
	(\Delta+k^2) G(x,k) = \delta(x), \ \ x\in\mathbb{R}^d, \ \ k \in \mathbb{C}^d\setminus \mathbb{R}^d;
\end{equation}
formula (\ref{4.1}) at fixed $k$ is considered as an equation for 
\begin{equation}
	\psi = e^{ikx}\mu(x,k),
\end{equation}
where $\mu$ is sought in $\mathbb{L}^{\infty}(\mathbb{R}^d)$; 
as a corollary of (\ref{4.1}), (\ref{4.2}), (\ref{4.8}), 
$\psi$ satisfies (\ref{eq}) for $E=k^2$; 
 $h$ of \eqref{4.4} is a generalized "`scattering"' amplitude. 

 
 In addition,  $h$, $\psi$, $G$  in their zero energy restriction, that is for $E=0$, were considered for the first time in \cite{Beals1985}.
The Faddeev functions $h$, $\psi$, $G$ were, actually, rediscovered in \cite{Beals1985}.

Let
\begin{equation}
\begin{aligned}
		\Sigma_E = \left\{ k\in \mathbb{C}^d: k^2 = k_1^2 + \ldots + k_d ^2 = E\right\},\\
		\Theta_E = \left\{ k\in \Sigma_E,\  l\in\Sigma_E: \mbox{Im}\, k = \mbox{Im}\, l\right\},\\
		|k| = (|\mbox{Re}\,k|^2 +|\mbox{Im}\,k|^2)^{1/2}.
\end{aligned}		
\end{equation}
Let  
\begin{equation}
	\begin{aligned}
	\text{$v$ satisfy \eqref{assumption},} \ \ \  \|v\|_{m,1}  \leq N,
	\end{aligned}
\end{equation}
		\begin{equation}
			\hat{v}(p) = (2\pi)^{-d}\int\limits_{\mathbb{R}^d} e^{ipx} v(x)dx, \ \ p\in \mathbb{R}^d,
		\end{equation}
then we have that:
\begin{equation}\label{mu_1}
	\mu(x,k) \rightarrow 1 \ \  \text{ as } \ \ |k|\rightarrow \infty 
\end{equation}
and, for any $\sigma>1$,
\begin{equation}\label{mu_2}
	|\mu(x,k)| + |\nabla\mu(x,k)|  \leq \sigma \ \  \text{ for } \ \ |k| \geq \lambda_1(N,m,d,r,\sigma),
\end{equation}
where $x\in \mathbb{R}^d$, $k \in \Sigma_E$;
\begin{equation}\label{lim_1}
			\hat{v}(p) = \lim\limits_
		{\scriptsize
			\begin{array}{c}
			(k,l)\in \Theta_E,\, k-l=p\\
			|\mbox{Im}\,k|=|\mbox{Im}\,l|\rightarrow \infty
			\end{array}
		} h(k,l)\ \ \ 	 \text{ for any } p\in \mathbb{R}^d,
	\end{equation}
\begin{equation}\label{lim_2}
	\begin{aligned}
		|\hat{v}(p) - h(k,l)|\leq \frac{c_1(m,d,r)N^2}{(E+\rho^2)^{1/2}}  	\ \ \text{ for }  (k,l) \in \Theta_E, \ \  p = k-l,\\
		|\mbox{Im}\,k| = |\mbox{Im} \,l| = \rho, \ \ \ E+\rho^2 \geq \lambda_2(N,m,d,r), \\ p^2 \leq 4(E+\rho^2).
	\end{aligned}
\end{equation}

Results of the type (\ref{mu_1}), (\ref{mu_2}) go back to \cite{Beals1985}.  
For more information concerning (\ref{mu_2}) see estimate (4.11) of \cite{IN2012}.
Results of the type (\ref{lim_1}), (\ref{lim_2}) 
(with less precise right-hand side in (\ref{lim_2})) go back to \cite{Henkin1987}. 
Estimate (\ref{lim_2}) follows, for example, from 
formulas (\ref{4.1}), (\ref{4.4}) and 
the estimate
\begin{equation}\label{4.15.new}
	\begin{aligned}
	\| \Lambda^{-s} g(k) \Lambda^{-s}\|_{\mathbb{L}^2(\mathbb{R}^d)\rightarrow \mathbb{L}^2(\mathbb{R}^d)} = O(|k|^{-1}) \\ \text{ as } \ |k|\rightarrow \infty,\ \ \
		k\in \mathbb{C}^d\setminus \mathbb{R}^d,
	\end{aligned}
\end{equation}
for $s>1/2$, where $g(k)$ denotes the integral operator with the Schwartz kernel $g(x-y,k)$ and $\Lambda$ denotes the multiplication operator by the function $(1+|x|^2)^{1/2}$. Estimate (\ref{4.15.new}) was formulated, first, in \cite{LN1987} for $d\geq 3$. Concerning proof of (\ref{4.15.new}), see \cite{Weder1991}. 

In addition, we have that:
\begin{equation}\label{delta_h0}
	\begin{aligned}
	h_2(k,l) - h_1(k,l) = (2\pi)^{-d} \int\limits_{\mathbb{R}^d} \psi_1(x,-l) (v_2(x) - v_1(x)) \psi_2(x,k) dx 	\\ 	 
	\text{ for }  (k,l) \in \Theta_E,\ |\mbox{Im}\,k| = |\mbox{Im} \,l| \neq 0,\\ \text{ and $v_1$, $v_2$ satisfying (\ref{4.6}),}
	\end{aligned}
\end{equation}
and, under assumtions of Theorem \ref{Theorem_2.1},
\begin{equation}\label{delta_v}
	\begin{aligned}
		|\hat{v}_1(p) - \hat{v}_2(p) - h_1(k,l) + h_2(k,l)|
		\leq \frac{
		c_2(m,d,r)N
		\|v_1-v_2\|_{\mathbb{L}^{\infty}(B_r)}
		}
		{(E+\rho^2)^{1/2}}  	\\ \text{ for }  (k,l) \in \Theta_E, \ \  p = k-l,\ \
		|\mbox{Im}\,k| = |\mbox{Im} \,l| = \rho,\\ E+\rho^2 \geq \lambda_3(N,m,d,r), \ \ p^2 \leq 4(E+\rho^2),
	\end{aligned}
\end{equation}
where $h_j$, $\psi_j$ denote $h$ and $\psi$ of (\ref{4.4}) and (\ref{4.1}) for $v = v_j$,  $j=1,2$.

Formula (\ref{delta_h0}) was given in \cite{Novikov1996}.
Estimate \eqref{delta_v} was given e.g. in \cite{IN2012++}.
\section{Proof of Theorem \ref{Theorem_2.1}}

Let
\begin{equation}\label{6.3}
\begin{aligned}
	\mathbb{L}^{\infty}_{\mu}(\mathbb{R}^d) = \{u\in \mathbb{L}^{\infty}(\mathbb{R}^d): \|u\|_\mu<+\infty\},\\
	\|u\|_\mu = \mbox{ess}\,\sup\limits_{p\in\mathbb{R}^d} (1+|p|)^{\mu}|u(p)|, \ \ \ \mu>0.
\end{aligned}
\end{equation}

Note that
\begin{equation}\label{6.4}
	\begin{aligned}
	w \in \mathbb{W}^{m,1}(\mathbb{R}^d) \Longrightarrow \hat{w} \in 	\mathbb{L}^{\infty}_{\mu}(\mathbb{R}^d)\cap { C}(\mathbb{R}^d),\\
		\|\hat{w}\|_\mu \leq c_3(m,d) \|w\|_{m,1} \ \ \ \text{ for } \ \ \mu = m, 
	\end{aligned}
\end{equation}
where $\mathbb{W}^{m,1}$,  $\mathbb{L}^{\infty}_{\mu}$ are the spaces of (\ref{2.5}), (\ref{6.3}),
\begin{equation}
	\hat{w}(p) = (2\pi)^{-d}\int\limits_{\mathbb{R}^d} e^{ipx}w(x) dx, \ \ \ p\in \mathbb{R}^d.
\end{equation}

Using the inverse Fourier transform formula
\begin{equation}
	w(x) = \int\limits_{\mathbb{R}^d} e^{-ipx}\hat{w}(p) dp, \ \ \ x\in \mathbb{R}^d,
\end{equation}
we have that
\begin{equation}\label{5.4}
	\begin{aligned}
		\|v_1 - v_2\|_{\mathbb{L}^{\infty}(B_r)} \leq 
		\sup\limits_{x\in \overline{B}_r}|\int\limits_{\mathbb{R}^d} e^{-ipx}\left(\hat{v}_2(p) - \hat{v}_1(p)\right) dp| \leq\\
		\leq I_1(\kappa) + I_2(\kappa) \ \ \ \text{ for any }  \ \kappa>0,
	\end{aligned}
\end{equation}
where
\begin{equation}\label{6.7}
	\begin{aligned}
		I_1(\kappa) = \int\limits_{|p|\leq \kappa} |\hat{v}_2(p) - \hat{v}_1(p)| dp, \\
		I_2(\kappa) = \int\limits_{|p|\geq \kappa} |\hat{v}_2(p) - \hat{v}_1(p)| dp.
	\end{aligned}
\end{equation}

Using (\ref{6.4}), we obtain that
\begin{equation}\label{6.8}
	|\hat{v}_2(p) - \hat{v}_1(p)| \leq 2c_3(m,d) N (1+|p|)^{-m}, \ \ \ p\in \mathbb{R}^d.
\end{equation}
Let
\begin{equation}
	c_4 =  \int\limits_{p\in \mathbb{R}^d,  |p|= 1} d p.
\end{equation}
Combining (\ref{6.7}), (\ref{6.8}), we find that, for any $\kappa>0$,
\begin{equation}\label{5.16}
	\begin{aligned}
		I_2(\kappa) \leq 2 c_3(m,d) N c_4 \int\limits_{\kappa}\limits^{+\infty} \frac{dt}{t^{m-d+1}} \leq \frac{2c_3(m,d)Nc_4}{m-d}\frac{1}{\kappa^{m-d}}.
	\end{aligned}
\end{equation}

Due to (\ref{delta_v}), we have that
\begin{equation}\label{6.9}
\begin{aligned}
	|\hat{v}_2(p) - \hat{v}_1(p)| \leq |h_2(k,l) - h_1(k,l)| + 
	\frac{c_2(m,d,r)N\|v_1-v_2\|_{\mathbb{L}^{\infty}(B_r)}}{(E+\rho^2)^{1/2}},\\
	\text{ for }  (k,l) \in \Theta_E, \ \  p = k-l,\ \
		|\mbox{Im}\,k| = |\mbox{Im} \,l| = \rho,\\ E+\rho^2 \geq \lambda_3(N,m,d,r), \ \ p^2 \leq 4(E+\rho^2).
\end{aligned}
\end{equation}

Let 
\begin{equation}
	\begin{aligned}
	\delta = 	||S_1^+(E) - S_2^+(E)||_{\mathbb{L}^{2}(\partial B_r \times \partial B_r)}.
	\end{aligned}
\end{equation}
Combining   \eqref{eq4.2}, \eqref{Id} and (\ref{delta_h0}), we get that
\begin{equation}\label{6.11}
	\begin{aligned}
	|h_2(k,l) - h_1(k,l)| \leq  \delta 
	\left\|\frac{\partial{\phi_1}}{\partial \nu_+} - 
	\frac{\partial{\phi_1}}{\partial \nu_-}\right\|_{\mathbb{L}^{2}(B_r)} 
	\left\|\frac{\partial{\phi_2}}{\partial \nu_+} - 
	\frac{\partial{\phi_2}}{\partial \nu_-}\right\|_{\mathbb{L}^{2}(B_r)},\\
	(k,l)\in \Theta_E, \ |\mbox{Im}\,k| = |\mbox{Im}\,l| \neq 0,
	\end{aligned}
\end{equation}
where $\phi_j$, $j=1,2$, denotes the solution of \eqref{4.3} with $v = v_j$, satisfying  
\begin{equation}
 \phi_j(x) = \psi_j(x,k) \ \ \ \text{ for } \ \ \ x\in  \overline{B}_{r}.  
\end{equation}
Using \eqref{eq4.3}, (\ref{mu_2})  and the fact that $C^1(\partial B_r) \subset \mathbb{H}^1(\partial B_r)$, we find that
\begin{equation}\label{6.12}
\begin{aligned}
	\left\|\frac{\partial{\phi_j}}{\partial \nu_+} - 
	\frac{\partial{\phi_j}}{\partial \nu_-}\right\|_{\mathbb{L}^{2}(B_r)} \leq
	 \sigma c_5(r,d)(1+E)\,\exp\bigg(|\mbox{Im}\,k|(r+1)\bigg), \\
	k \in \Sigma_E, \ |k|\geq \lambda_1(N,m,d,r,\sigma), \ j=1,2.
\end{aligned}
\end{equation}
Here and bellow in this section the constant $\sigma$ is the same that in (\ref{mu_2}).
 
Combining (\ref{6.11}) and (\ref{6.12}), we obtain that
\begin{equation}\label{6.13}
	\begin{aligned}
		|h_2(k,l) - h_1(k,l)| \leq c_5^2 \sigma^2 (1+E)^{2}  e^{2\rho (r+1)} \delta, \\ \text{ for }  
		(k,l)\in\Theta_E, \ \ \rho = |\mbox{Im}\,k|= |\mbox{Im}\,l|, \\ E+\rho^2\geq \lambda_1^2(N,m,d,r,\sigma).
	\end{aligned}
\end{equation}
Using (\ref{6.9}), (\ref{6.13}), we get that
\begin{equation}\label{6.14}
	\begin{aligned}
		|\hat{v}_2(p) - \hat{v}_1(p)| \leq c_5^2 \sigma^2 (1+E)^{2} e^{2\rho (r+1)} \delta +\\
		+ \frac{c_2(m,d,r)N\|v_1-v_2\|_{\mathbb{L}^{\infty}(B_1)}}{(E+\rho^2)^{1/2}},\\
		p\in\mathbb{R}^d,\  p^2 \leq 4(E+\rho^2),\  E+\rho^2 \geq \max\{ \lambda_1^2 , \lambda_3\}.
	\end{aligned}
\end{equation}
Let 
\begin{equation}
	\varepsilon = \left(\frac{1}{2c_2(m,d,r)N c_6}\right)^{1/d}, \ \ \ c_6 = \int\limits_{p\in\mathbb{R}^d, |p|\leq 1} d p,
\end{equation}
and $\lambda_4(N,m,d,r,\sigma)>0$ be such that
\begin{equation}\label{5.14new.}	
	E+\rho^2 \geq \lambda_4(N,m,d,r,\sigma) \Longrightarrow 
	\left\{
	\begin{aligned}
	&E+\rho^2 \geq \lambda_1^2(N,m,d,r,\sigma), \\ &E+\rho^2 \geq \lambda_3(N,m,d,r), \\ &
	\left(\varepsilon(E+\rho^2)^{\frac{1}{2d}}\right)^2 \leq 4(E+\rho^2).
	\end{aligned}\right.
\end{equation} 
Using (\ref{6.7}), (\ref{6.14}), we get that
\begin{equation}\label{5.15}
\begin{aligned}
	I_1(\kappa) \leq c_6 \kappa^d \Big( c_5^2 \sigma^2 (1+E)^{2} e^{2\rho (r+1)} \delta
	+ \frac{c_2(m,d,r)N\|v_1-v_2\|_{\mathbb{L}^{\infty}(B_1)}}{(E+\rho^2)^{1/2}}\Big),\\  \kappa>0, \ \kappa^2\leq 4(E+\rho^2), \\ 
	E+\rho^2 \geq \lambda_4(N,m,d,r,\sigma).
\end{aligned}
\end{equation}
 Combining (\ref{5.4}),  (\ref{5.16}), (\ref{5.15}) for $\kappa = \varepsilon(E+\rho^2)^{\frac{1}{2d}}$ and (\ref{5.14new.}), we get that
 \begin{equation}\label{5.18}
 \begin{aligned}
 	\|v_1 - v_2\|_{\mathbb{L}^{\infty}(B_r)} \leq c_7(N,m,d,r,\sigma) (1+E)^{2}\sqrt{E+\rho^2}\, e^{2\rho (r+1)} \delta+ \\+ 
 	c_8(N,m,d)(E+\rho^2)^{-\frac{m-d}{2d}} + \frac{1}{2}\|v_1 - v_2\|_{\mathbb{L}^{\infty}(B_r)}, \\
 	E+\rho^2 \geq \lambda_4(N,m,d,r,\sigma).
 	\end{aligned}
 \end{equation}
 
Let $\tau' \in (0,1)$,
\begin{equation}\label{5.19}
	\beta = \frac{1-\tau'}{2(r+1)}, \ \ \ \rho = \beta\ln\left(3+ \delta^{-1}\right),
\end{equation} 
and $\delta_1 = \delta_1(N,m,d,\sigma,r,\tau')>0$ be such that 
\begin{equation}\label{new24}	
	\delta \in (0,\delta_1)\Longrightarrow 
	\left\{
	\begin{aligned}
	&E+\left(\beta\ln\left(3+ \delta^{-1}\right)\right)^2 \geq \lambda_4(N,m,d,r,\sigma), \\ 
	&E+\left(\beta\ln\left(3+ \delta^{-1}\right)\right)^2 \leq 	(1+E)\left(\beta\ln\left(3+ \delta^{-1}\right)\right)^2,	\\ 
	\end{aligned}\right.
\end{equation}
 Then for the case when $\delta \in (0,\delta_1)$, due to (\ref{5.18}), we have that
\begin{equation}\label{5.20}
	\begin{aligned}
		\frac{1}{2}\|v_1 - v_2\|_{\mathbb{L}^{\infty}(B_r)} &\leq \\ \leq 
		c_7 (1+E)^2 &\left(E+ \left(\beta\ln\left(3+ \delta^{-1}\right)\right)^2\right)^{\frac{1}{2}}
		\left(3+ \delta^{-1}\right)^{2\beta (r+1)} \delta + \\
		&+  c_8 \left(E+ \left(\beta\ln\left(3+ \delta^{-1}\right)\right)^2\right)^{-\frac{m-d}{2d}} =\\
		= c_7 (1+E)^2  &\left(E+ \left(\beta\ln\left(3+ \delta^{-1}\right)\right)^2\right)^{\frac{1}{2}}
		\left(1+ 3\delta \right)^{1-\tau'} \delta^{\tau'} 
		 + \\
		&+ c_8 \left(E+ \left(\beta\ln\left(3+ \delta^{-1}\right)\right)^2\right)^{-\frac{m-d}{2d}}.
	\end{aligned}
\end{equation}
Combining  \eqref{new24} and \eqref{5.20}, we obtain that  for  $s \in [0, s^*]$, $\tau \in (0,\tau')$ and $\delta \in (0,\delta_1)$
the following estimate holds:
\begin{equation}\label{new25}
	\begin{aligned}
||v_2 - v_1||_{L^\infty(B_r)} \leq c_{9}(1+E)^{\frac{5}{2}} \delta^\tau + 
	c_{10} (1+E)^{\frac{s- s^*}{2}}\left(\ln\left(3+\delta ^{-1}\right)\right)^{-s}, 
	\end{aligned}
\end{equation}
where $s^* = \frac{m-d}{d}$ and  $c_{9}, c_{10} >0$ depend only on  $N$,   $m$, $d$,  $r$, $\sigma$, $\tau'$ and $\tau$.

Estimate (\ref{new25}) in the general case (with modified $c_{9}$ and $c_{10}$) follows from
 (\ref{new25})  for $\delta \leq \delta_1(N,m,d,\sigma,r,\tau')$ and and the property that 
\begin{equation}\label{new23}
\|v_j\|_{\mathbb{L}^{\infty}(B_r)} \leq c_{11}(m,d)N.
\end{equation}

This completes the proof of \eqref{eq_t2} 


\section{Buckhgeim-type analogs of the Faddeev functions} 
Let us identify $\mathbb{R}^2$ with $\mathbb{C}$ and use coordinates
$z = x_1+ix_2$,  $\bar{z} = x_1-ix_2$, where $(x_1,x_2) \in \mathbb{R}^2$.
Following \cite{NS2010}-\cite{S2011}, we consider the functions $G_{z_0}, \psi_{z_0}, \tilde{\psi}_{z_0}$, $\delta h_{z_0}$ 
going back to Buckhgeim's paper \cite{Buckhgeim2008} and being analogs of the Faddeev functions:
\begin{equation}\label{eq_bu_psi}
	\begin{aligned}
	\psi_{z_0}(z,\lambda) = e^{\lambda(z-z_0)^2} + \int\limits_{B_{r}} G_{z_0}(z,\zeta,\lambda)
	v(\zeta) \psi_{z_0}(\zeta,\lambda)\, d\mbox{Re}\zeta\, d\mbox{Im}\zeta,
	\\
	\widetilde{\psi}_{z_0}(z,\lambda) = e^{\bar{\lambda}(\bar{z}-\bar{z}_0)^2} + \int\limits_{B_{r}} \overline{G_{z_0}(z,\zeta,\lambda)}
	v(\zeta) \widetilde{\psi}_{z_0}(\zeta,\lambda)\, d\mbox{Re}\zeta\, d\mbox{Im}\zeta,
	\end{aligned}
\end{equation}
\begin{equation}
	\begin{aligned}\label{Bu_G}
		G_{z_0}(z,\zeta,\lambda) = \frac{1}{4\pi^2}\int\limits_{B_{r}} 
		\frac
		{e^{-\lambda(\eta-z_0)^2+\bar{\lambda}(\bar{\eta}-\bar{z_0})^2}d\mbox{Re}\eta\,d\mbox{Im}\eta}
		{(z-\eta)(\bar{\eta}-\bar{\zeta})} \, e^{\lambda(z-z_0)^2-\bar{\lambda}(\bar{\zeta}-\bar{z_0})^2},\\
		z = x_1 + ix_2,\  z_0\in B_{r},\  \lambda \in \mathbb{C},
	\end{aligned}
\end{equation} 
where  $v$ satisfies (\ref{assumption2});
\begin{equation}\label{deltah}
	\delta h_{z_0}(\lambda) = \int\limits_{B_{r}} \widetilde{\psi}_{z_0,1}(z,-\lambda) 
	\left(v_2(z) - v_1(z)\right)
	\psi_{z_0,2}(z,\lambda) \, d\mbox{Re}z\, d\mbox{Im}z, \ \ \lambda \in \mathbb{C},
\end{equation}
where $v_1$, $v_2$ satisfy (\ref{assumption2}) and $\widetilde{\psi}_{z_0,1}$, $\psi_{z_0,2}$ denote 
$\widetilde{\psi}_{z_0}$, $\psi_{z_0}$ of (\ref{eq_bu_psi}) for $v=v_1$ and $v=v_2$, respectively.

We recall that (see \cite{NS2010}, \cite{NS2011}):
\begin{itemize}

\item The function $G_{z_0}$ satisfies the equations
\begin{equation}\label{5.4b}
	\begin{aligned}
	4\frac{\partial^2}{\partial z \partial \bar{z}} G_{z_0}(z,\zeta,\lambda) = \delta(z - \zeta),\\
	4\frac{\partial^2}{\partial \zeta \partial \bar{\zeta}} G_{z_0}(z,\zeta,\lambda) = \delta(z - \zeta),
	\end{aligned}
\end{equation}
where $z,z_0, \zeta \in B_{r}$, $\lambda \in \mathbb{C}$ and $\delta$ is the Dirac delta function; 

\item Formulas (\ref{eq_bu_psi}) 
at fixed $z_0$ and $\lambda$ are
considered as equations for $\psi_{z_0}$, $\widetilde{\psi}_{z_0}$ in $L^{\infty}(B_{r})$; 

\item As a corollary of (\ref{eq_bu_psi}),
(\ref{Bu_G}), (\ref{5.4b}), the functions $\psi_{z_0}$, $\widetilde{\psi}_{z_0}$ satisfy (\ref{eq}) in $B_{r}$ for $E=0$ and $d=2$;

\item The function
$\delta h_{z_0}$ is similar to  the right side of (\ref{delta_h0}). 
\end{itemize}

Let potentials $v,v_1,v_2 \in C^2(\overline{B}_{r})$ and 
\begin{equation}
	\begin{aligned}
	\|v\|_{C^2(\overline{B}_{r})}  \leq N, \ \ \|v_j\|_{C^2(\overline{B}_{r})}  \leq N, \ \ j=1,2,\\
	(v_1-v_2)|_{\partial B_{r}} = 0, \ \  \ \frac{\partial }{\partial \nu}(v_1-v_2)|_{\partial B_{r}} = 0,
	\end{aligned}
\end{equation}
then we have that:
\begin{equation}\label{5.5b}
	\begin{aligned}
		\psi_{z_0}(z,\lambda) = e^{\lambda(z-z_0)^2} \mu_{z_0}(z,\lambda), \ \ \
		\widetilde{\psi}_{z_0}(z,\lambda) = e^{\bar{\lambda}(\bar{z}-\bar{z}_0)^2} \widetilde{\mu}_{z_0}(z,\lambda),
	\end{aligned}
\end{equation} 
\begin{equation}\label{5.6b}
	\mu_{z_0}(z,\lambda) \rightarrow 1, \ \ \widetilde{\mu}_{z_0}(z,\lambda) \rightarrow 1 \ \ \text{ as } |\lambda| \rightarrow \infty
\end{equation}
and, for any $\sigma>1$,
\begin{subequations} \label{5.7b}
	\begin{equation}
		|\mu_{z_0}(z,\lambda)| + |\nabla\mu_{z_0}(z,\lambda)| \leq \sigma, 
	\end{equation}
	\begin{equation}
	|\widetilde{\mu}_{z_0}(z,\lambda)| + |\nabla\widetilde{\mu}_{z_0}(z,\lambda)|  \leq \sigma,
	\end{equation}
\end{subequations} 
where $\nabla = \left(\partial/\partial x_1, \partial/\partial x_2\right)$, $z= x_1 + ix_2$, $z_0 \in B_{r}$, $\lambda \in \mathbb{C}$, $|\lambda| \geq \rho_1(N,r,\sigma)$;
\begin{equation}\label{5.8b}
	\begin{aligned}
	v_2(z_0) - v_1(z_0) = \lim\limits_{\lambda\rightarrow \infty} \frac{2}{\pi}|\lambda| \delta h_{z_0}(\lambda)\\
		\text{ for any } z_0\in B_{r},
	\end{aligned}
\end{equation}
\begin{equation}\label{5.9b}
\begin{aligned}
	\left| v_2(z_0) - v_1(z_0) -  \frac{2}{\pi}|\lambda| \delta h_{z_0}(\lambda) \right| \leq 
	\frac{c_{12}(N,r) \left(\ln(3|\lambda|)\right)^2}
	{|\lambda|^{3/4}}\\
		\text{ for } z_0\in B_{r}, \ |\lambda| \geq \rho_2(N,r). 
\end{aligned}
\end{equation}
Formulas (\ref{5.5b}) can be considered as definitions of $\mu_{z_0}$, $\widetilde{\mu}_{z_0}$. Formulas (\ref{5.6b}), (\ref{5.8b})
were given in \cite{NS2010}, \cite{NS2011} and go back to \cite{Buckhgeim2008}. Estimates \eqref{5.7b} were proved in \cite{IN2012}.
Estimate (\ref{5.9b}) was obtained in \cite{NS2010}, \cite{S2011}. 
\section{Proof of Theorem \ref{Theorem_2.2}}

We suppose that 
 $\widetilde{\psi}_{z_0,1}(\cdot,-\lambda)$, $\psi_{z_0,2}(\cdot,\lambda)$, $\delta h_{z_0}(\lambda)$
are defined as in Section 6 but with $v_j-E$ in place of $v_j$, $j=1,2$. 
Note that functions $\widetilde{\psi}_{z_0,1}(\cdot,-\lambda)$, $\psi_{z_0,2}(\cdot,\lambda)$ 
satisfy (\ref{eq}) in $B_r$ with $v = v_j$, $j=1,2$, respectively.
We also use the notation $N_E = N+E$. Then, using (\ref{5.9b}), we have that
\begin{equation}\label{7.1}
\begin{aligned}
	\left| v_2(z_0) - v_1(z_0) -  \frac{2}{\pi}|\lambda| \delta h_{z_0}(\lambda) \right| \leq 
	\frac{c_{12}(N_E,r) \left(\ln(3|\lambda|)\right)^2}
	{|\lambda|^{3/4}}\\
		\text{ for } z_0\in B_r, \ |\lambda| \geq \rho_2(N_E,r). 
\end{aligned}
\end{equation}
Let 
\begin{equation}
	\delta  = ||S_1^+(E) - S_2^+(E)||_{\mathbb{L}^{2}(\partial B_r \times \partial B_r)}.
\end{equation}
Combining  \eqref{eq4.2}, \eqref{Id} and (\ref{deltah}),  we get that
\begin{equation}\label{new7.3}
	\begin{aligned}
	|\delta h_{z_0}(\lambda) | \leq  \delta 
	\left\|\frac{\partial{\phi_1}}{\partial \nu_+} - 
	\frac{\partial{\phi_1}}{\partial \nu_-}\right\|_{\mathbb{L}^{2}(B_r)} 
	\left\|\frac{\partial{\phi_2}}{\partial \nu_+} - 
	\frac{\partial{\phi_2}}{\partial \nu_-}\right\|_{\mathbb{L}^{2}(B_r)},\\
	(k,l)\in \Theta_E, \ |\mbox{Im}\,k| = |\mbox{Im}\,l| \neq 0,
	\end{aligned}
\end{equation}
where $\phi_j$, $j=1,2$, denotes the solution of \eqref{4.3} with $v = v_j$, satisfying  
\begin{equation}
	\begin{aligned}
 \phi_1(x) = \widetilde{\psi}_{z_0,1}(x,-\lambda),\ \ \
 \phi_2(x) = \psi_{z_0,2}(x,\lambda),
  \ \ \ \text{ for }  x\in  \overline{B}_{r}.
 \end{aligned}  
\end{equation}
Using \eqref{eq4.3},  (\ref{5.7b}) and the fact that $C^1(\partial B_r) \subset \mathbb{H}^1(\partial B_r)$, we find that:
\begin{equation}\label{7.5}
\begin{aligned}
	\left\|\frac{\partial{\phi_j}}{\partial \nu_+} - 
	\frac{\partial{\phi_j}}{\partial \nu_-}\right\|_{\mathbb{L}^{2}(B_r)} \leq
	 \sigma c_{13}(r)(1+E)\,\exp\bigg(|\lambda|(4r^2+4r)\bigg), \\ 
		\lambda \in \mathbb{C}, \ \  |\lambda|\geq\rho_1(N_E,r,\sigma), \ \ j=1,2.
\end{aligned}
\end{equation}
Here and bellow in this section the constant $\sigma$ is the same that in (\ref{5.7b}).

Combining (\ref{new7.3}), (\ref{7.5}), we obtain that
\begin{equation}\label{7.6}
	\begin{aligned}
	|\delta h_{z_0}(\lambda)| \leq  c_{14}(E,r,\sigma) \,\exp\bigg(|\lambda|(8r^2+8r)\bigg) \delta, \\
		\lambda \in \mathbb{C}, \ \  |\lambda|\geq\rho_1(N_E,r,\sigma).
	\end{aligned}
\end{equation}
Using (\ref{7.1}) and (\ref{7.6}), we get that
\begin{equation}\label{7.7}
	\begin{aligned}
		\left| v_2(z_0) - v_1(z_0) \right| \leq  c_{14}(E,r,\sigma) \,\exp\bigg(|\lambda|(8r^2+8r)\bigg) \delta
		+ \\+	\frac{c_{12}(N_E,r) \left(\ln(3|\lambda|)\right)^2}
	{|\lambda|^{3/4}}, \\
		z_0 \in B_r, \ \ \lambda \in \mathbb{C}, \ \  |\lambda|\geq\rho_3(N_E,r,\sigma) = \max\{\rho_1,\rho_2\}.
	\end{aligned}
\end{equation}

We fix some $\tau \in (0,1)$ and let
\begin{equation}\label{7.8}
	\beta = \frac{1-\tau}{8r^2+8r}, \ \ \ \lambda = \beta\ln\left(3+ \delta^{-1}\right),
\end{equation} 
where $\delta$ is so small that $|\lambda| \geq \rho_3(N_E,r,\sigma)$. Then due to (\ref{7.7}), we have that
\begin{equation}\label{7.9}
	\begin{aligned}
		\|v_1 - v_2\|_{\mathbb{L}^{\infty}(B_r)} &\leq 
	c_{14}(E,r,\sigma) 
		\left(3+ \delta^{-1}\right)^{\beta(8r^2+8r)} \delta + \\
	&+  c_{12}(N_E,r)  	\frac{
	\left(\ln\left(3\beta\ln\left(3+ \delta^{-1}\right)\right)\right)^2}
	{\left(\beta\ln\left(3+ \delta^{-1}\right)\right)^{\frac{3}{4}} }  =\\
		&= c_{14}(E,r,\sigma)  \left(1+ 3\delta \right)^{1-\tau} \delta^{\tau} 
		 + \\
		&+ c_{12}(N_E,r)  \beta^{-\frac{3}{4}}	\frac{
	\left(\ln\left(3\beta\ln\left(3+ \delta^{-1}\right)\right)\right)^2}
	{\left(\ln\left(3+ \delta^{-1}\right)\right)^{\frac{3}{4}} },
	\end{aligned}
\end{equation}
where  $\tau, \beta$ and $\delta$ are the same as in (\ref{7.8}). 

Using (\ref{7.9}), we obtain that
\begin{equation}\label{7.10}
	\|v_1 - v_2\|_{\mathbb{L}^{\infty}(B_r)} \leq c_{15}(N,E,r,\sigma)\left(\ln\left(3+\delta ^{-1}\right)\right)^{-\frac{3}{4}} 
	\left(\ln\left(3\ln\left(3+ \delta^{-1}\right)\right)\right)^2
\end{equation}
for $\delta = ||S_1^+(E) - S_2^+(E)||_{\mathbb{L}^{2}(\partial B_r \times \partial B_r)} \leq \delta_2(N_E,r,\sigma)$, where
$\delta_2$ is a sufficiently small positive constant. 
Estimate (\ref{7.10})  in the general case (with modified $c_{15}$) follows from 
(\ref{7.10}) for $\delta \leq \delta_2(N_E,r,\sigma)$ 
and the property that $\|v_j\|_{\mathbb{L}^{\infty}(B_r)} \leq N$.

This completes the proof of (\ref{estimation_2}).


\section{Proof of Lemma \ref{Lemma_4.1}}
In this section we assume for simplicity that $r = 1$ and therefore $\partial B_r = \mathbb{S}^{d-1}$.

We fix an orthonormal basis in $\mathbb{L}^2(\partial B_r)$:
\begin{equation}
	\begin{array}{l}
	\displaystyle
		\{f_{jp} : j \geq 0;\  1 \leq p \leq p_j \}, \\
		\text{$f_{jp}$ is a spherical harmonic of degree $j$,}
	\end{array}
\end{equation}
   where $p_j$  is the dimension of the space of spherical harmonics of order $j$,
\begin{equation}
	p_j = \binom {j+d-1} {d-1} - \binom {j+d-3} {d-1},
\end{equation}
	where 
\begin{equation}
	\binom {n} {k} = 
		\frac{n(n-1)\cdots(n-k+1)}{k!} \ \ \ \text{ for $n \geq 0$} 
\end{equation}	
and
\begin{equation}
	\binom {n} {k} = 
		0 \ \ \ \text{ for $n < 0$.} 
\end{equation}
The precise choice of $f_{jp}$ is irrelevant for our purposes. Besides orthonormality, we only need
$f_{jp}$ to be the restriction of a homogeneous harmonic polynomial of degree $j$ to the sphere $\partial B_r$
and so $|x|^j f_{jp}(x/|x|)$ is harmonic pn $\mathbb{R}^d$. In the Sobolev spaces $\mathbb{H}^s(\partial B_r)$ the norm is defined by
\begin{equation}\label{lp5}
	\left\| \sum_{j,p} c_{jp}f_{jp}\right\|^2_{\mathbb{H}^s(\partial B_r)} = \sum_{j,p}(1+j)^{2s}|c_{jp}|^2.
\end{equation} 
The solution $\phi$ of the exterior Dirichlet problem
\begin{equation}\label{lp6}
	\begin{aligned}
	-\Delta\phi    = E\phi, \ \ \text{ in } \mathbb{R}^d \setminus \overline{B}_r, \\
	\lim\limits_{|x|\rightarrow +\infty} |x|^{\frac{d-1}{2}} \left(\frac{\partial}{\partial|x|}\phi - i\sqrt{E}\phi\right)=0,\\
	\phi|_{\partial B_r} = u \in \mathbb{H}^1(\partial B_r),
	\end{aligned}
\end{equation}
can be expressed in the following form (see, for example, \cite{Berezanskii1958}, \cite{CK1998}):
\begin{equation}\label{lp7}
	\phi = \sum_{j,p} c_{jp}\phi_{jp},
\end{equation}
where $c_{jp}$ are expansion coefficients of  $u$ in the basis $\{f_{jp}: j \geq 0;\  1 \leq p \leq p_j\}$,
 and
\begin{equation}\label{lp8}
	\begin{aligned}
		&\text{$\phi_{jp}$ denotes the solution of \eqref{lp6} with $u=f_{jp}$},\\
	&\phi_{jp}(x) = h_{jp}(|x|)f_{jp}(x/|x|),\\
	&h_{jp}(|x|) = |x|^{-\frac{d-2}{2}}
	\frac{H^{(1)}_{j+\frac{d-2}{2}}(\sqrt{E}|x|)} 
	{H^{(1)}_{j+\frac{d-2}{2}}(\sqrt{E})},
	\end{aligned}
\end{equation}
where $H_\mu^{(1)}$ is the Hankel function of the first kind.
Let
\begin{equation}\label{lp9}
	\phi^0_{jp}(x) = |x|^{-j-d+2}f_{jp}(x/|x|).
\end{equation}
Note that $\phi^0_{jp}$ is harmonic in $\mathbb{R}^d \setminus \{0\}$ and 
\begin{equation}\label{lp10}
	\lim\limits_{|x|\rightarrow +\infty} |x|^{\frac{d-1}{2}} \left(\frac{\partial}{\partial|x|}\phi^0_{jp} - i\sqrt{E}\phi^0_{jp}\right)=0
	\ \ \ \text{ for  $j+\frac{d-3}{2}>0$}.
\end{equation}
Using the Green formula and the radiation condition for $\phi_{jp}$, $\phi^0_{jp}$, we get that 
\begin{equation}\label{lp11}
	\begin{aligned}
 \int\limits_{\mathbb{R}^d\setminus B_r} E \phi_{jp} \phi_{jp}^0 dx =  
 \int\limits_{\mathbb{R}^d\setminus B_r}  \left( \Delta \phi_{jp}^0 \phi_{jp} - \Delta \phi_{jp} \phi_{jp}^0  \right) dx =\\
 =\int\limits_{\partial B_r} 
 \left(\frac{\partial\phi_{jp}^0}{\partial \nu_+} \phi_{jp} - \frac{\partial\phi_{jp}}{\partial \nu_+} \phi_{jp}^0\right) dx 	\ \ \ \text{ for  $j+\frac{d-3}{2}>0$}.
 \end{aligned}
\end{equation}
Due to  \eqref{lp8} and \eqref{lp9}, we have that
\begin{equation}\label{lp12}
	\left|\int\limits_{\partial B_r} 
	\frac{\partial\phi_{jp}^0}{\partial \nu_+} \phi_{jp} dx\right| = (j+d - 2) \int\limits_{\partial B_r} f_{jp}^2 dx = j+d - 2.
\end{equation}
Using also the following property of the Hankel function of the first kind (see, for example, \cite{Watson1944}):
\begin{equation}\label{lp13}
	|H^{(1)}_\mu(x)| \text{ is a decreasing function of $x$ for $x\in \mathbb{R}_+, \mu \in \mathbb{R}$,}
\end{equation} 
we get that 
\begin{equation}\label{lp14}
	\begin{aligned}
	\left|\int\limits_{\mathbb{R}^d\setminus B_r}  \phi_{jp} \phi_{jp}^0 dx\right| = 
	\left|\int\limits_{1}\limits^{+\infty} t^{-j-d+2} h_{jp}(t) t^{d-1} dt\right| =\\ =
	\left|\int\limits_{1}\limits^{+\infty} t^{-j-\frac{d}{2}} \frac{H^{(1)}_{j+\frac{d-2}{2}}(\sqrt{E}t)} 
	{H^{(1)}_{j+\frac{d-2}{2}}(\sqrt{E})}dt\right| \leq \int\limits_{1}\limits^{+\infty} t^{-j-\frac{d}{2}} dt
	= \frac{1}{j+\frac{d}{2}-1} \leq 2 \\  \text{ for  $j+\frac{d-3}{2} >0$}.
	\end{aligned}
\end{equation}
Combining \eqref{lp8}, \eqref{lp9},  \eqref{lp11}, \eqref{lp12} and \eqref{lp14}, we obtain that
\begin{equation}\label{lp15}
	\left|\int\limits_{\partial B_r}  \frac{\partial\phi_{jp}^0}{\partial \nu_+} \phi_{jp} dx \right| =
	\left|\frac{h_{jp}'(r)}{h_{jp}(r)}\right| \leq j+d - 2 + 2E \ \ \ \text{ for  $j+\frac{d-3}{2}>0$}.
\end{equation}

Let consider the cases when $j+\frac{d-3}{2} \leq 0$.

\noindent
{\it Case 1.} $j=0$, $d=2$. Using the property $d H^{(1)}_0(t) /dt = - H^{(1)}_1(t)$,  we get that
\begin{equation}\label{lp16}
	\frac{h_{jp}'(r)}{h_{jp}(r)} = \sqrt{E}\frac{H^{(1)}_1(\sqrt{E})}{H^{(1)}_0(\sqrt{E})}.
\end{equation}
We recall that functions $H^{(1)}_0$ and $H^{(1)}_1$ have the following asymptotic forms (see, for example \cite{Watson1944}):
\begin{equation}\label{lp17}
	\begin{aligned}
		H^{(1)}_0(t) \sim \frac{2i}{\pi} \ln (t/2) \ \ \text{ as } \ \ t\rightarrow +0,\\
		H^{(1)}_0(t) \sim \sqrt{\frac{2}{\pi t}} e^{i(t-\pi/4)} \ \ \text{ as } \ \ t\rightarrow +\infty,\\
		H^{(1)}_1(t) \sim -\frac{i}{\pi}  (2/t) \ \ \text{ as } \ \ t\rightarrow +0,\\
		H^{(1)}_1(t) \sim \sqrt{\frac{2}{\pi t}} e^{i(t-3\pi/4)} \ \ \text{ as } \ \ t\rightarrow +\infty.
	\end{aligned}
\end{equation}
Using \eqref{lp13} and \eqref{lp17}, we get that for some $c>0$
\begin{equation}\label{lp18}
	\frac{H^{(1)}_1(t)}{H^{(1)}_0(t)} \leq c\left(1 + 1/t\right).
\end{equation}
Combining \eqref{lp16} and \eqref{lp18}, we obtain that for $j=0$, $d=2$
\begin{equation}\label{lp19}
	\left|\frac{h_{jp}'(r)}{h_{jp}(r)}\right| \leq c(1 + \sqrt{E}).
\end{equation}

\noindent
{\it Case 2.} $j=0$, $d=3$. We have that
\begin{equation}\label{lp20}
	H_{j+\frac{d-2}{2}}^{(1)}(t) = \sqrt{\frac{2}{\pi t}} e^{i(t-\pi/2)}. 
\end{equation}
Using \eqref{lp8} and \eqref{lp20}, we get that for $j=0$, $d=3$
\begin{equation}\label{lp21}
	\frac{h_{jp}'(r)}{h_{jp}(r)} = -1 + i\sqrt{E}.
\end{equation}

Combining  \eqref{lp5}-\eqref{lp8}, \eqref{lp15}, \eqref{lp19} and \eqref{lp21}, we get that for some constant $c'=c'(d)>0$
	\begin{equation}\label{lp22}
		\left\|\frac{\partial \phi}{\partial \nu_+}\right\|_{\mathbb{L}^2(\partial B_r)}^2 = 
		\sum_{j,p} c_{jp}^2 \left|\frac{ h_{jp}'(r)}{h_{jp}(r)}\right|^2 \leq c'(1+E)^2\sum_{j,p} (1+j)^2 		c_{jp}^2.
	\end{equation}
	Using \eqref{lp5} and \eqref{lp22}, we obtain \eqref{eq4.3}
\section*{Acknowledgements}
This work was fulfilled in the framework of research carried out under the supervision of R.G. Novikov.
This work was partially supported by FCP Kadry No. 14.A18.21.0866.








\end{document}